\documentclass[12pt]{article}
\usepackage[a4paper, total={7in, 10in}]{geometry}
\usepackage[parfill]{parskip}
\usepackage{physics, tensor, float, subcaption}
\usepackage{graphicx}
\graphicspath{ {Plots/} }
\usepackage{jhep-mod}
\usepackage{bm}
\usepackage{soul}
\usepackage{amssymb,amsmath,amsthm}
\usepackage{mathrsfs}
\usepackage[utf8]{inputenc}
\usepackage{enumerate}
\usepackage{bigints}
\usepackage{xcolor}
\usepackage{appendix}
\usepackage{graphicx}
\usepackage{float}
\usepackage{tikz}
\usepackage{setspace}
\usepackage{cancel}
\usepackage{doi}
\definecolor{purple}{rgb}{1,0,1}
\definecolor{lime}{HTML}{A6CE39} 

\definecolor{lime}{HTML}{A6CE39}
\newcommand{\orcidicon}{%
	\begin{tikzpicture}
	\draw[lime, fill=lime] (0,0) 
		circle [radius=0.16] 
		node[white] {{\fontfamily{qag}\selectfont \tiny ID}};
	\draw[white, fill=white] (-0.0625,0.095) 
		circle [radius=0.007];
	\end{tikzpicture}
	\hspace{-5mm}
}
\newcommand\orcidMatt{{\href{https://orcid.org/0000-0003-1088-6485}{\orcidicon}}}

\renewcommand{\O}{\mathcal{O}}
\begin{document}

\title{\null\vspace{-50pt}
\leftline{Primes between 8th powers?}
}

\author{
\Large
Matt Visser\!\orcidMatt\!
}
\affiliation{School of Mathematics and Statistics, Victoria University of Wellington, \\
\null\qquad PO Box 600, Wellington 6140, New Zealand.}
\emailAdd{matt.visser@sms.vuw.ac.nz}
\def\theta{\vartheta}
\def\O{{\mathcal{O}}}
\def\Li{{\mathrm{Li}}}

\abstract{
\vspace{1em}

OpenAI (using a non-public ``internal frontier model'') has recently (6 October 2026) made a number of remarkable claims, including some claims regarding the quasi-Riemann hypothesis (QRH). 
While there are good reasons to be cautious regarding these claims, one can nevertheless ask what (if true) their implications would be. 
While there is widespread agreement that any version of the QRH will provide various improved bounds on the location of the prime numbers, I feel that there is some value in making these improvements explicit. 
Specifically, if we assume validity of the claimed 7/8 version of the QRH then I shall demonstrate (for all $n\in \mathbb{Z^+}$):\\
1) the existence of primes between eighth powers, in the interval $(n^8,(n+1)^8)$;\\
2)  the existence of primes in the shorter interval $(n^8,n^8+8n^7)$;\\
3)  the existence of primes in the even shorter interval $(n^8,n^8+n^7)$.\\
Similarly for reals greater than unity $x>1$ one can demonstrate: \\
4) the existence of primes in the Hoheisel-style interval $(x,x+x^{7/8})$.\\
All of this is contingent on the validity of the claimed 7/8 version of the QRH.

\bigskip

\bigskip
\noindent
{\sc Date:} 14 October 2026; \LaTeX-ed \today

}

\maketitle
\def\tr{{\mathrm{tr}}}
\def\diag{{\mathrm{diag}}}
\def\cof{{\mathrm{cof}}}
\def\pdet{{\mathrm{pdet}}}
\def\d{{\mathrm{d}}}
\parindent0pt
\parskip7pt
\def\Kerr{{\scriptscriptstyle{\mathrm{Kerr}}}}
\def\eos{{\scriptscriptstyle{\mathrm{eos}}}}
\section{Introduction}

OpenAI (using a non-public ``internal frontier model'') has recently uploaded some 722 preprints onto {\sf github}~\cite{slopocalypse}. There are good reasons to be cautious regarding the remarkable claims being made:
\begin{itemize}
\item Most of the 722 preprints are unreadable and would earn an automatic desk reject if submitted to a scientific journal.
\item At least 3 of the preprints have already been withdrawn due to computational errors~\cite{withdrawn}, and 12 other preprints have been corrected~\cite{withdrawn}.
\item Some 42\% of the 722 preprints are claimed to have been formalized in the Lean language, and then axiomatically verified for formal correctness. The  issue here is that there is significant evidence that the initial autoformalization process, that is, automatically translating natural language (English) to formal statements (Lean), is not always reliable~\cite{autoformalization}.
\item
As of 14 October 2026  trust in generic AI-generated ``proofs'' is not helped by other AIs continuing to spout manifest nonsense about the Riemann hypothesis implying the Legendre hypothesis.
\item
As of 14 October 2026  trust in at least some (allegedly) human commentators  on these issues is not helped by at least some of these commentators not even realizing that the zeros of the Riemann zeta function are symmetrically distributed around the line $\Re(\rho)=1/2$. 
\end{itemize}

\clearpage
Specifically, OpenAI makes two claims  regarding the quasi-Riemann hypothesis:
\begin{description}
\item[QRH$_{11/12}$:] \ \\
That the zeros $\rho$ of the Riemann zeta function satisfy $\Re(\rho)\leq 11/12$ \cite{QRH1}.
(Though not stated this way, this would imply $\Re(\rho) \in [1/12,11/12]$.)
\item[QRH$_{7/8}$: ] \ \\
That the zeros $\rho$ of the Riemann zeta function satisfy $\Re(\rho)\leq 7/8$ \cite{QRH2}.\\
(Though not stated this way, this would imply $\Re(\rho) \in [1/8,7/8]$.)
\end{description}
There have already been two announced updates on these topics:
\begin{description}
\item[QRH$_{47/48}$:] 
A potentially simplified proof~\cite{simplified} that the zeros $\rho$ of the Riemann zeta function satisfy $\Re(\rho)\leq 47/48$.
(Though not stated this way, this would imply $\Re(\rho) \in [1/48,47/48]$.)
\item[QRH$_{34999/40000}$:] 
A (very) slightly strengthened result~\cite{strengthened}, asserting that the zeros $\rho$ of the Riemann zeta function satisfy $\Re(\rho)\leq {34999\over 40000}= {7/8}- {1\over 40000}$.
(Though not stated this way, this would imply $\Re(\rho) \in [{5001\over40000}, {34999\over 40000}]\subset [1/8,7/8]$.)
\end{description}
The situation is very rapidly evolving, sometimes on a daily basis.

In the current document I shall focus on the {\bf QRH}$_{7/8}$ version of these claims.\\
Assuming validity of {\bf QRH}$_{7/8}$,  I shall demonstrate that for all $n\in \mathbb{Z^+}$:\\
1) there are primes between eighth powers, in the interval $(n^8,(n+1)^8)$;\\
2)  there are primes in the shorter interval $(n^8,n^8+8n^7)$;\\
3)  there are primes in the even shorter interval $(n^8,n^8+n^7)$.\\
Furthermore \\
4) for all real $x>1$ there are primes in the interval $(x,x+x^{7/8})$.

\section{Updating the Chebyshev functions}

In a very recent preprint \cite{Chebyshev2}, Fiori and Jaskari have shown, (assuming  the {\bf QRH}$_{7/8}$ version of the OpenAI claims is correct), that the second Chebyshev function $\psi(x)$ satisfies:
\begin{equation}
|\psi(x)-x| < 0.04934 \;x^{3/4} \;\ln(x)^2 + 0.02636425 \;x^{7/8}; \qquad (x\geq e^3).
\end{equation}
Sometimes it is convenient to slightly weaken the region of validity to $x > 21$.

For current purposes, it is more useful to work with the first Chebyshev function $\theta(x)$ and note~\cite{Chebyshev1}
\begin{equation}
0 < \psi(x)-\theta(x) < {\sqrt{x} \; \ln(x)^2\over 2 \ln 2}; \qquad (x>0).
\end{equation}
That is
\begin{equation}
\psi(x) -{\sqrt{x} \; \ln(x)^2\over 2 \ln 2}  <\theta(x)<\psi(x); \qquad (x>0).
\end{equation}
This allows us, for $x>e^3$,  to write the \emph{asymmetric} bounds
\begin{equation}
\label{E:theta_low}
\theta(x) > x -0.04934 \;x^{3/4} \ln(x)^2 - 0.02636425 \;x^{7/8} - {\sqrt{x} \;\ln(x)^2\over 2 \ln 2};
\end{equation}
and
\begin{equation}
\label{E:theta_high}
\theta(x) < x +0.04934 \;x^{3/4} \ln(x)^2 + 0.02636425 \; x^{7/8}.
\end{equation}

\section{Primes between eighth powers}

{\bf Theorem:}
Assuming {\bf QRH}$_{7/8}$,
 then for all positive integers $n$ we have
\begin{equation}
\pi((n+1)^8) > \pi(n^8); \qquad (n \in\mathbb{Z^+}).
\end{equation}

\begin{proof}
From equation \eqref{E:theta_low} we see
\begin{eqnarray}
\theta((x+1)^8) &> &
(x+1)^8  -3.15776 \; x^{6} \; \ln(x+1)^2 - 0.02636425 \; (x+1)^{7} 
\nonumber\\
&&\qquad - 32 {(x+1)^4 \; \ln(x+1)^2\over \ln 2}.
\end{eqnarray}
This holds for $(x+1)^8 \geq e^3$, that is for $x > e^{3/8}-1$, so certainly for $x\geq 1$.

From equation \eqref{E:theta_high} we see
\begin{equation}
\label{E:theta_high_8}
\theta(x^8) < 
x^8  +3.15776 \; x^{6} \; \ln(x)^2 + 0.02636425 \; x^{7}.
\end{equation}
This holds for $x^8 \geq e^3$, that is $x \geq e^{3/8}$, so certainly for $x\geq 2$.

Thence for $x\geq 2$ we have
\begin{eqnarray}
&&\theta((x+1)^8) - \theta(x^8) > (x+1)^8 - x^8
\nonumber\\
&&\qquad
-3.15776 \;(x+1)^{6} \; \ln(x+1)^2 - 0.02636425\; (x+1)^{7} - {32 \; (x+1)^4 \; \ln(x+1)^2\over \ln 2}
\nonumber\\
&&
\qquad -3.15776 \;x^{6} \ln(x)^2 - 0.02636425 \; x^{7}.
\end{eqnarray}
The RHS of this inequality is easily checked to be  positive for $x>2.473174837$, certainly for 
$x\geq 3$, whence
\begin{equation}
\theta((x+1)^8) - \theta(x^8) > 0; \qquad (x \geq 3);
\end{equation}
implying
\begin{equation}
\pi((x+1)^8) - \pi(x^8) > 0; \qquad (x \geq 3).
\end{equation}
Explicitly checking smaller integers
\begin{equation}
\pi(1)=0; \qquad \pi(2^8) = 54; \qquad \pi(3^8) = 847;
\end{equation}
we see that for all positive integers $n$ we have
\begin{equation}
\pi((n+1)^8) > \pi(n^8); \qquad (n \in\mathbb{Z^+}).
\end{equation}
So, assuming the validity of  {\bf QRH}$_{7/8}$, there are always primes between eighth powers.
\end{proof}
\section{Primes in the interval $[n^8, n^8+8 n^7]$}

A slightly stronger result can be obtained by considering the slightly narrower interval $[n^8, n^8+8 n^7] \subset [n^8,(n+1)^8]$. 

{\bf Theorem:}
Assuming {\bf QRH}$_{7/8}$,
 for all positive integers $n$ we have
\begin{equation}
\pi(n^8+8n^7) > \pi(n^8); \qquad (n \in\mathbb{Z^+}).
\end{equation}

\begin{proof}
From equation \eqref{E:theta_low} we see
\begin{eqnarray}
\theta(x^8+8x^7) &>& x^8+8x^7 -0.04934 \;(x^8+8x^7)^{3/4} \; \ln(x^8+8x^7)^2 - 0.02636425 \; (x^8+8x^7)^{7/8} 
\nonumber\\
&&\qquad - {\sqrt{x^8+8x^7} \; \ln(x^8+8x^7)^2\over 2 \ln 2}.
\end{eqnarray}
This holds for $x^8+8x^7> e^3$, 
that is $x>1.119407260$,
so certainly holds for $x\geq 2$.

We still have \eqref{E:theta_high_8} 
\begin{equation}
\theta(x^8) < 
x^8  +3.15776 x^{6} \ln(x)^2 + 0.02636425 x^{7}; \qquad (x\geq 2).
\end{equation}

Consequently for $x\geq 2$ we have
\begin{eqnarray}
&&\theta(x^8+8x^7) - \theta(x^8)  > 8x^7\nonumber\\
&&\qquad\qquad
-0.04934 \;(x^8+8x^7)^{3/4} \; \ln(x^8+8x^7)^2 - 0.02636425 \; (x^8+8x^7)^{7/8} 
\nonumber\\
&&\qquad\qquad
- {\sqrt{x^8+8x^7} \; \ln(x^8+8x^7)^2\over 2 \ln 2}
\nonumber\\
&&\qquad\qquad
-3.15776 \; x^{6} \; \ln(x)^2 - 0.02636425 \; x^{7}.
\end{eqnarray}
The RHS of this inequality is positive for 
$x > 4.459691950$, and so is certainly positive for 
$x\geq 5$.

Thence
\begin{equation}
\theta(x^8+8x^7) - \theta(x^8)  > 0; \qquad (x\geq 5);
\end{equation}
implying
\begin{equation}
\pi(x^8+8x^7) > \pi(x^8); \qquad (x\geq 5).
\end{equation}
It is then easy to check by explicit computation that this also holds for the larger range $x\geq 1$, and so 
\begin{equation}
\pi(x^8+8x^7) > \pi(x^8); \qquad (x\geq 1).
\end{equation}
In particular, for all positive integers we have 
\begin{equation}
\pi(n^8+8n^7) > \pi(n^8); \qquad (n\in\mathbb{Z^+}).
\end{equation}
\end{proof}

\section{Primes in the interval $[n^8, n^8+ n^7]$}

With a bit more effort a slightly stronger result can be obtained by considering the even narrower interval $[n^8, n^8+ n^7] \subset [n^8, n^8+8 n^7] \subset [n^8,(n+1)^8]$.

{\bf Theorem:}
Assuming {\bf QRH}$_{7/8}$,
 for all positive integers $n$ we have
\begin{equation}
\pi(n^8+n^7) > \pi(n^8); \qquad (n \in\mathbb{Z^+}).
\end{equation}

\begin{proof}
From equation \eqref{E:theta_low} we see
\begin{eqnarray}
\theta(x^8+x^7) &>& x^8+x^7 -0.04934 \;(x^8+x^7)^{3/4} \;\ln(x^8+x^7)^2 - 0.02636425 \;(x^8+x^7)^{7/8} 
\nonumber\\
&&\qquad - {\sqrt{x^8+x^7} \;\ln(x^8+x^7)^2\over 2 \ln 2}.
\end{eqnarray}
This holds for $x^8+x^7\geq e^3$, that is $ x > 1.358017028$, so this certainly holds for $x\geq 2$.

We still have \eqref{E:theta_high_8} 
\begin{equation}
\theta(x^8) < 
x^8  +3.15776 \; x^{6} \; \ln(x)^2 + 0.02636425 \; x^{7}; \qquad (x\geq 2).
\end{equation}

Consequently for $x\geq 2$ we have
\begin{eqnarray}
&&\theta(x^8+x^7) - \theta(x^8)  > x^7\nonumber\\
&&\qquad\qquad
-0.04934 \; (x^8+x^7)^{3/4} \; \ln(x^8+x^7)^2 - 0.02636425\;(x^8+x^7)^{7/8} 
\nonumber\\
&&\qquad\qquad
- {\sqrt{x^8+x^7} \;\ln(x^8+x^7)^2\over 2 \ln 2}
\nonumber\\
&&\qquad\qquad
-3.15776 \; x^{6} \; \ln(x)^2 - 0.02636425 \;x^{7}.
\end{eqnarray}
The RHS of this inequality is positive for $x >180.4094594$, 
so it is certainly positive for $x\geq 181$.

That is
\begin{equation}
\theta(x^8+x^7) - \theta(x^8)  > 0; \qquad (x\geq 181);
\end{equation}
implying
\begin{equation}
\pi(x^8+x^7) - \pi(x^8)  > 0; \qquad (x\geq 181);
\end{equation}

Checking the range $[1,181]$ is now somewhat trickier.\\
Let us immediately specialise to positive integers $n$.
\begin{itemize}
\item That $\pi(n^8+n^7) - \pi(n^8)>0$ for positive integers $n\in[1,16]$ can easily be checked by direct computation on a modern laptop.
\item
For integers in the range $n\in[16,181]$ it is more practical to reformulate the problem and check the inequalities
\begin{equation}
\hbox{nextprime}(n^8) < n^8+ n^7.
\end{equation}
(Calculations of $\hbox{nextprime}(n^8)$ are computationally much easier than directly evaluating the prime counting function.)
\end{itemize}
Performing these checks, for all positive integers we have 
\begin{equation}
\pi(n^8+n^7) > \pi(n^8); \qquad (n\in\mathbb{Z^+}).
\end{equation}
\end{proof}

\section{Primes in the Hoheisel-style interval $[x, x+ x^{7/8}]$}

With yet more work we can look for primes in the 
Hoheisel-style interval~\cite{Hoheisel,Visser} given by $[x, x+ x^{7/8}]$.
This is largely (but not entirely) a repackaging of our previous result for intervals of the form $[x^8, x^8+x^7]$. Some delicacy is required when explicitly checking ``low'' numbers.

{\bf Theorem:}
Assuming {\bf QRH}$_{7/8}$, for positive reals $x$ we have
\begin{equation}
\pi(x+ x^{7/8}) > \pi(x); \qquad (x > 1).
\end{equation}

\begin{proof}
From equation \eqref{E:theta_low} we see
\begin{eqnarray}
\theta(x+ x^{7/8}) &>& x+ x^{7/8} -0.04934 \;(x+ x^{7/8})^{3/4} \;\ln(x+ x^{7/8})^2 - 0.02636425 \; (x+ x^{7/8})^{7/8} 
\nonumber\\
&&\qquad - {\sqrt{x+ x^{7/8}} \; \ln(x+ x^{7/8})^2\over 2 \ln 2}.
\end{eqnarray}
This holds for $x+ x^{7/8}> e^3$, that is $ x > 11.56755903$, so this certainly holds for $x\geq 12$.

We still have \eqref{E:theta_high} 
\begin{equation}
\theta(x) < 
x +0.04934 \; x^{3/4} \; \ln(x)^2 + 0.02636425 \; x^{7/8};
\qquad (x > e^3).
\end{equation}

Consequently for $x> e^3$, certainly for $x \geq 21$,  we have 
\begin{eqnarray}
&&\theta(x+x^{7/8}) - \theta(x)  > x^{7/8}\nonumber\\
&&\qquad\qquad
-0.04934 \;(x+x^{7/8})^{3/4} \;\ln(x+x^{7/8})^2 - 0.02636425 \;(x+x^{7/8})^{7/8} 
\nonumber\\
&&\qquad\qquad
- {\sqrt{x+x^{7/8}} \;\ln(x+x^{7/8})^2\over 2 \ln 2}
\nonumber\\
&&\qquad\qquad
-0.04934 \;x^{3/4} \ln(x)^2 - 0.02636425 \; x^{7/8}.
\end{eqnarray}
The RHS of this inequality is positive for $x >1.122210792\times10^{18} = (180.4094594)^8$, 
so it is certainly positive for $x\geq (9/8)\times10^{18}$.

That is
\begin{equation}
\theta(x+x^{7/8}) - \theta(x)  > 0; 
\qquad \left(x\geq (9/8)\times10^{18}\right);
\end{equation}
implying
\begin{equation}
\pi(x+x^{7/8}) - \pi(x)  > 0; 
\qquad \left(x\geq (9/8)\times10^{18}\right).
\end{equation}

Checking the range $[1,(9/8)\times10^{18}]$ is now very much trickier.\\
\begin{itemize}
\item That $\pi(x+x^{7/8}) - \pi(x)>0$ for positive reals $x\in[1,10^8]$ can easily be checked by direct computation on a modern laptop.
\item
For reals in the range $x\in [10^8, (9/8)\times10^{18}]$ we 
proceed as follows: \\
Note that $[x,x+x^{7/8}] \supset [x,x+\sqrt{x}]$,
and that $[x,x+\sqrt{x}]$ is the interval used in checking the upper half of Oppermann's conjecture, which can be written as
$\pi(x+\sqrt{x})>\pi(x)$. 
Opperman's conjecture has now been explicitly verified for $x<(3.33\times 10^{13})^2 = 1.10 \times 10^{27}$~\cite{Oppermann1}, which is already (much) more than enough for current purposes. (The published version~\cite{Oppermann2} expands the verification to $x<(7\times 10^{13})^2 = 4.9 \times 10^{27}$.) 
\end{itemize}
As a result of these checks, for all positive reals we have 
\begin{equation}
\pi(x+x^{7/8}) > \pi(x); \qquad (x>1).
\end{equation}
\end{proof}

\clearpage
\section{Conclusions}

OpenAI's claimed 7/8 version of the QRH~\cite{QRH2} would, if verified, have rather strong implications for the distribution of primes. Herein I have tried to make various versions of these implications as clear and explicit as possible. If one does \emph{not} use something like {\bf QRH}$_{7/8}$ then the best that can currently be said are much weaker statements  along the lines that:\\
 1) for $n\geq 1$ there will be primes between all subsequent 86-th powers~\cite{Lee}, \\
 see also references~\cite{Cully-Hugill,CHJ1,CHJ2};\\
 2)  for $n > \exp(\exp(33))$, there will be primes between all subsequent 3-rd powers (cubes)~\cite{Dudek,Cully-Hugill}, \\
 3) for unspecified sufficiently large $x$ there will be primes in the Hoheisel-style interval $[x,x+x^{0.525}]$, see~\cite{BHP}.

OpenAI's claimed results regarding the QRH are certainly remarkable~\cite{QRH1,QRH2}, with significant implications regarding the distribution of the primes.
This topic is likely to generate continuing ongoing interest.

\clearpage

\bigskip
\hrule\hrule\hrule

\begin{thebibliography}{99}
\newcommand{\arXiv}[1]{arXiv:~{\href{https://arxiv.org/abs/#1}{\color{blue}#1}}}

\bibitem{slopocalypse}
OpenAI, 
``Sharing AI progress in mathematics'',
\url{https://openai.com/index/sharing-ai-progress-in-mathematics/}\\
Specifically, see: 
\url{https://github.com/openai/math/tree/main/preprints},\\
(6 October 2026).

\bibitem{withdrawn}
\url{https://github.com/openai/math/blob/main/history.md}\\
(Three articles withdrawn, another twelve corrected, as of 14 October 2026.)

\bibitem{autoformalization}
Alexander Bastounis, Fabian Circelli, Anders C. Hansen,\\
``Navier-Stokes lost in translation: Why Lean verification of AI autoformalisation does not guarantee correct natural language proofs'', \\
\arXiv{2610.08144} [math.AP], (6 October 2026).


\bibitem{QRH1}
OpenAI, ``The quasi-Riemann hypothesis'',
\url{https://github.com/openai/math/blob/main/preprints/The-Quasi-Riemann-Hypothesis-October-5-2026/paper2.pdf},
(5 October 2026).


\bibitem{QRH2}
OpenAI, ``The quasi-Riemann hypothesis: A Zero-Free Half-Plane $\Re(s)>7/8$'',
 \url{https://github.com/openai/math/blob/main/preprints/The-Quasi-Riemann-Hypothesis-September-30-2026/paper.pdf},\\
(30 September 2026).


\bibitem{simplified}
Shiva Kintali,
``A Short Proof of the Quasi-Riemann Hypothesis'',
\url{https://shivakintali.github.io/papers/QRH.pdf},
(7 October 2026).


\bibitem{strengthened}
Baiying Liu,\\
``Slightly improved zero-free half-planes for the quasi-Riemann hypothesis'',
\arXiv{2610.12234} [math.NT], (8 October 2026).


\bibitem{Chebyshev2}
Andrew Fiori and Mikko Jaskari,
``Explicit bounds for the prime number theorem'',
\arXiv{2609.23222} [math.NT], (8 October 2026).

\bibitem{Chebyshev1}
See for example \url{https://en.wikipedia.org/wiki/Chebyshev_function}

\bibitem{Hoheisel}
Guido Hoheisel, 
``Primzahlprobleme in der Analysis'',\\
Sitzungsberichte Berliner Akad. d. Wiss., (1930) 580-588.

\bibitem{Visser}
Matt Visser, ``Effective short intervals containing primes'',\\
\arXiv{2508.18786} math.NT 



\bibitem{Oppermann1}
Jonathan Sorenson, Jonathan Webster,\\
``An algorithm and computation to verify Legendre's Conjecture up to $3.33\times 10^{13}$'',\\
\arXiv{2401.13753} [math.NT] (24 Jan 2024).

\bibitem{Oppermann2}
Jonathan Sorenson, Jonathan Webster,\\
``An algorithm and computation to verify Legendre's Conjecture up to $7\times 10^{13}$'',\\
 Research in Number Theory {\bf 11} (2025) 4,
\url{https://doi.org/10.1007/s40993-024-00589-4}

\clearpage
\bibitem{Lee}
Ethan Simpson Lee,
``Minimal zero-free regions for results on primes
between consecutive perfect k-th powers'',
\arXiv{2602.14340} [math.NT] (3 March 2026).\\
Acta Arith. 2026 (in press).

\bibitem{Cully-Hugill}
Michaela Cully-Hugill,
``Primes between consecutive powers'',\\
Journal of Number Theory, {\bf 247} (2023) 100--117.
\doi{10.1016/j.jnt.2022.12.002}
[\arXiv{2107.14468} [math.NT]]

\bibitem{CHJ1}
Michaela Cully-Hugill and Daniel R. Johnston,\\
``On the error term in the explicit formula of Riemann--von Mangoldt'',\\
 International Journal of Number Theory
{\bf 19 \# 6} (2023) 1205--1228
\doi{10.1142/S1793042123500598}
[\arXiv{2111.10001} [math.NT]]

\bibitem{CHJ2}
Michaela Cully-Hugill and Daniel R. Johnston,\\
``On the error term in the explicit formula of Riemann--von Mangoldt. II.'', \\
Funct. Approx. Comment. Math. 
 {\bf 73 \# 2} (2025) 223--242 . 
\doi{10.7169/facm/241110-18-11}
[\arXiv{2402.04272} [math.NT]]

\bibitem{Dudek}
Adrian Dudek,
``An Explicit Result for Primes Between Cubes'',\\
 Funct. Approx. Comment. Math. {\bf 55 \# 2} (2016) 177--197. 
\doi{10.7169/facm/2016.55.2.3}
[\arXiv{1401.4233} [math.NT]]

\bibitem{BHP}
R. C. Baker, G. Harman, and J. Pintz,\\
``The Difference Between Consecutive Primes, II",\\
Proc. London Math. Soc. {\bf 83 \# 3} (2001) 532--562.

\end{thebibliography}
\end{document}